\documentclass [12pt]{amsart}

\usepackage{amsmath, accents} %

\usepackage{amssymb,amsxtra,amsfonts}
\usepackage{amsmath, accents} %
\usepackage{graphics}

\usepackage[colorlinks]{hyperref}

\usepackage{scalerel}   %
\def\apple{{\includegraphics[height=2ex]{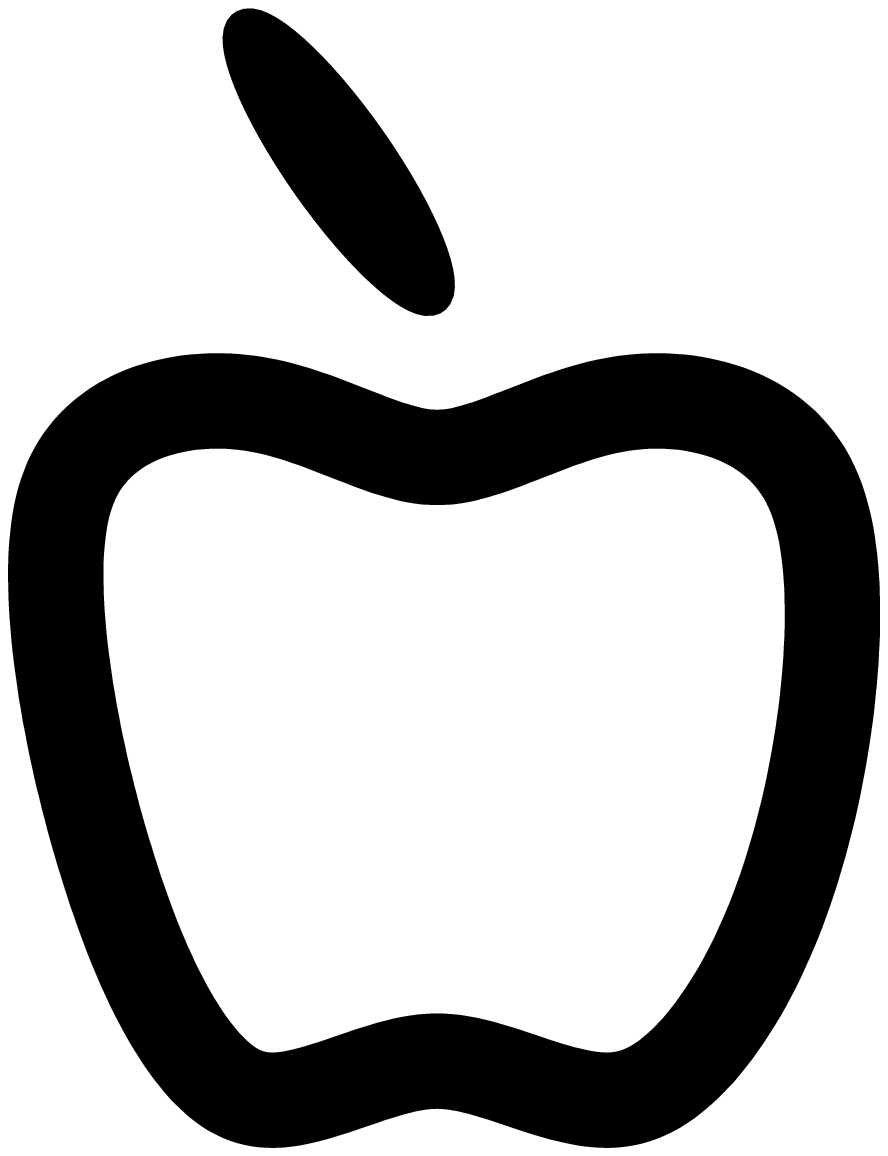}}}

\long\def\pb #1*/{}

\def\reE@DeclareMathSymbol#1#2#3#4{%
    \let#1=\undefined
    \DeclareMathSymbol{#1}{#2}{#3}{#4}}
\DeclareSymbolFont{symbolsC}{U}{txsyc}{m}{n}
\SetSymbolFont{symbolsC}{bold}{U}{txsyc}{bx}{n}
\DeclareFontSubstitution{U}{txsyc}{m}{n}
\reE@DeclareMathSymbol{\strictiff}{\mathrel}{symbolsC}{76}

\newcommand\beq{\begin{equation}}
\newcommand\eeq{\end{equation}}
\newcommand\bal{\begin{align*}}
\newcommand\eal{\end{align*}}   
\newcommand\bmx{\left(\begin{matrix}}
\newcommand\emx{\end{matrix}\right)}
\newcommand\bsmx{\left(\begin{smallmatrix}}
\newcommand\esmx{\end{smallmatrix}\right)}
\newcommand\bmxnp{\begin{matrix}}
\newcommand\emxnp{\end{matrix}}
\newcommand\bsmxnp{\begin{smallmatrix}}
\newcommand\esmxnp{\end{smallmatrix}}

\newcommand{\wt}{\widetilde}

\DeclareMathSymbol{\widehatsym}{\mathord}{largesymbols}{"62}
\newcommand\lowerwidehatsym{%
  \text{\smash{\raisebox{-1.3ex}{%
    $\widehatsym$}}}}
\newcommand\fixwidehat[1]{%
  \mathchoice
    {\accentset{\displaystyle\lowerwidehatsym}{#1}}
    {\accentset{\textstyle\lowerwidehatsym}{#1}}
    {\accentset{\scriptstyle\lowerwidehatsym}{#1}}
    {\accentset{\scriptscriptstyle\lowerwidehatsym}{#1}}
}

\newcommand{\wh}{\fixwidehat}

\newcommand{\bSi}{{\bf \Si}}

\newcommand{\monesp}{\;\!\!}   

\newcommand{\spq}{/\!\!/}
\newcommand{\sqp}{\setminus\monesp\monesp\!\!\setminus\,}

\providecommand{\deg}{\text{\rm deg}}

\providecommand{\Irr}{\text{\rm Irr}}

\providecommand{\<}{\langle}
\renewcommand{\>}{\rangle}

\def\part#1{\frac{\partial\phantom{q}}{\partial#1}}

\newcommand{\MDR}{\mathcal{M}_{\text{\rm DR}}}

\newcommand{\MB}{\mathcal{M}_{\text{\rm B}}}

\newcommand{\MDol}{\mathcal{M}_{\text{\rm Dol}}}

\newcommand{\Id}{\text{\rm Id}}

\DeclareMathOperator{\ISto}{{\IS}to} 
 
\DeclareMathOperator{\rk}{\mathop{\rm rk}}
\newcommand{\ram}{\mathop{\rm ram}}

\newcommand{\Prod}{\prod}

\DeclareMathOperator{\Hom}{Hom}         

\newcommand{\SL}{{\mathop{\rm SL}}}

\newcommand{\GL}{{\mathop{\rm GL}}}

\DeclareMathOperator{\End}{End}

\newcommand{\bd}{{\bf d}}

\DeclareSymbolFont{bbold}{U}{bbold}{m}{n}
\DeclareSymbolFontAlphabet{\mathbbold}{bbold}

\newcommand{\IA}{\mathbb{A}}

\newcommand{\IC}{\mathbb{C}}

\newcommand{\IN}{\mathbb{N}}

\newcommand{\IP}{\mathbb{P}}

\newcommand{\IS}{\mathbb{S}}

\newcommand{\IZ}{\mathbb{Z}}

\newcommand{\cA}{\mathcal{A}}
\newcommand{\cB}{\mathcal{B}}

\newcommand{\cC}{\mathcal{C}}

\newcommand{\cE}{\mathcal{E}}

\newcommand{\cI}{\mathcal{I}}

\newcommand{\cM}{\mathcal{M}}

\newcommand{\gM}{       \mathfrak{M}     }

\renewcommand{\sl}{       \mathfrak{sl}     } 

\newcommand{\be}{\beta}

\newcommand {\eps}{\varepsilon}

\newcommand{\Ga}{\Gamma}

\newcommand{\Si}{\Sigma}
\newcommand{\Th}{\Theta}

\makeatletter
 \newlength{\typesize}
\makeatother

\newlength{\vvoff}
\newlength{\hhoff}

\newcommand{\pf}{\begin{bpf}}

\newcommand{\pfms}{\begin{bpfms}}
\newcommand{\epf}{\end{bpf}\hfill$\square$\\}           
\newcommand{\epfms}{\end{bpfms}\hfill$\square$\\}       

\newcommand{\idea}{\begin{bidea}}

\newcommand{\eidea}{\end{bidea}\hfill$\square$\\}           

\newcommand{\sk}{\begin{bsk}}    

\newcommand{\esk}{\end{bsk}\hfill$\square$\\}           
\newcommand{\sketch}{\begin{bsketch}}

\newcommand{\esketch}{\end{bsketch}\hfill$\square$\\}

\newtheorem {hypo}{\bf\hspace{-\parindent}Hypothesis}
\newtheorem {thm}[hypo]{Theorem}   
\newtheorem {prop}[hypo]{Proposition}

\newtheorem {lem}[hypo]{Lemma}

\newtheorem {defn}[hypo]{Definition}

\theoremstyle{remark}\newtheorem{rmk}[hypo]{Remark}
\theoremstyle{remark}

\begin{document}


\title{Diagrams for nonabelian Hodge spaces on the affine line}
\author{Philip Boalch and Daisuke Yamakawa}%

\begin{abstract}
In this announcement a diagram will be defined for each nonabelian Hodge space on the affine line. 
\end{abstract}

\maketitle

\section{Introduction}

In a previous paper (\cite{rsode} Apx. C), a diagram  was defined for each algebraic connection on a vector bundle on the affine line, under the condition that the connection was untwisted at infinity (i.e. had unramified irregular type).
In that case the diagram was a graph (or a ``doubled quiver'').
It is known that such a connection determines a triple of complex algebraic moduli spaces
$\MDR, \MDol, \MB$, which are algebraically distinct yet naturally diffeomorphic. In other words the connection determines
a single {\em nonabelian Hodge space} $\gM,$ with a triple of distinct algebraic structures.
Hence this gives a way to attach a diagram to a class of 
nonabelian Hodge spaces.
The simplest examples of such graphs match up with the 
affine Dynkin diagrams corresponding to the Okamoto Weyl group symmetries of the Painlev\'e equations (corresponding to some of the H3 surfaces, i.e. the spaces $\gM$ of real dimension four).
More background and motivation related to integrable systems and isomonodromy is recounted in \cite{hit70}.

The purpose of this note is to extend this story by defining 
a diagram  
for any algebraic connection on a vector bundle on the affine line, i.e. for any nonabelian Hodge space attached to the affine line.
One can show using the Fourier--Laplace transform 
that any moduli space $\MDR$ of meromorphic connections on a smooth affine curve of genus zero is isomorphic to one on the affine line 
(i.e. with just one puncture), and this is expected to hold for the full nonabelian Hodge triple. It is thus hoped that these diagrams serve a useful purpose in the classification of 
nonabelian Hodge spaces (and this is certainly the case in the examples considered in \cite{rsode, slims, cmqv}).

\section{The construction}

Let $\Si=\IP^1(\IC)$ and fix a point $p=\infty\in \Si$ so that 
$\Si^\circ=\Si\setminus p$ is the affine line. 
The diagram of an algebraic connection $(E,\nabla)\to \Si^\circ$
is determined by its formal isomorphism class at $\infty$.
This formal class is equivalent to the irregular class plus the formal monodromy conjugacy classes, as follows.

Let 
$\partial$ be the circle of real oriented directions at 
$p$.
Recall that the exponential local system $\cI$ 
is a covering space $\pi:\cI\to \partial$, 
consisting of a disjoint union of circles $\<q\>$ each of which is a finite cover of $\partial$.
(This notation is from \cite{twcv, tops} where further 
discussion and references may be found).
Deligne's way of stating the Hukuhara--Turritin--Levelt 
formal classification of connections is as follows:

\begin{thm}[\cite{malg-book} Thm 2.3]
The category of connections on vector bundles on the formal punctured disk at $p$ is equivalent to the category of $\cI$-graded 
local systems $V^0\to \partial$ of finite dimensional complex vector spaces.
\end{thm}

Such a graded local system $V^0$ is the same thing as a 
local system (of finite dimensional complex  vector spaces) on the
topological space $\cI$, 
having {\em compact support} (in the sense that it has rank zero on all but a finite number of component circles of $\cI$).

By definition  the {\em irregular class} of a connection is the map $\Th:\pi_0(\cI) \to \IN$ taking the rank of $V^0$ on each circle 
(\cite{twcv} \S3.5).
In down to earth terms the circles $\<q\>$ correspond to the exponential factors $e^q$ of the corresponding connection, so fixing the irregular class amounts to fixing the exponential factors plus their integer multiplicities.
Thus an irregular class can be written as a formal sum
$$\Th = %
n_1\<q_1\>+\cdots+n_m\<q_m\>$$
of a finite number of distinct circles $\<q_i\>$, with integer multiplicities $n_i=\Th(q_i)\ge 1$.

Now a local system of rank $n$ on a circle $\<q\>$ is a very simple object, and is classified by the conjugacy class 
in $\GL_n(\IC)$ of its monodromy (in a positive sense once around the circle).
Thus the graded/formal  local system 
$V^0$ determines conjugacy classes 
$$\cC=(\cC_1,\ldots,\cC_m),\qquad 
\cC_i\subset \GL_{n_i}(\IC)$$
where $\cC_i$ is the class of the monodromy of the 
local system on the circle  $\<q_i\>$.

Now there is a well-known method 
due to Kraft--Procesi and others 
of attaching a graph $L$ (a type $A$ Dynkin graph)
to a {\em marked} conjugacy class in $\GL_n(\IC)$.
It is reviewed  in \cite{cmqv} Defn 9.2. 
Moreover the graph is independent of the marking if the marking
is chosen to be {\em minimal}, in the sense of \cite{cmqv} Defn 9.2.
The number of nodes of $L$ is the degree of a minimal polynomial of any element of the class. 
Thus $V^0$ determines {\em legs} $L_1,\ldots,L_m$, where
$L_i$ is the type $A$ Dynkin graph 
determined by a minimal marking of the class $\cC_i$.

The next step is to define the {\em core diagram} $\Ga$.
For this first recall (e.g. from \cite{tops}) that:

1) For any circle $\<q\>\subset \cI$ the ramification degree $\ram(q)$ is the degree of the covering 
map $\pi:\<q\>\to \partial$,

2) The set of {\em points of maximal decay} is a discrete subset  
$\apple{}\subset \cI$. It consists of  
a finite subset $\apple(q)\subset \<q\>$ in each circle,  where the function $e^q$ has maximal decay,

3) The size of the set $\apple(q)$ is called the irregularity $\Irr(q)$ of the irregular class $\<q\>$ (and is zero if and only if $q=0$).
For an arbitrary class $\Th=\sum n_i\<q_i\>$ the irregularity is
$\Irr(\Th)=\sum n_i\Irr(q_i)$,

4) For any pair of circles $\<q_1\>,\<q_2\>$ the irregular class 
$\Hom(\<q_1\>,\<q_2\>)$ is well-defined (the definition is straightforward if one thinks in terms of corresponding graded local systems).

Now the {\em core diagram} $\Ga$
is defined as follows:

$\bullet$ $\Ga$ has $m$ nodes,  
labelled by the circles $\<q_1\>,\ldots,\<q_m\>$,

$\bullet$ If $i\neq j$ then the number of arrows from 
$\<q_i\>$ to $\<q_j\>$ is given by 
\beq\label{eq: Bij} B_{ij}=A_{ij}-\be_i\be_j\eeq
where $\be_i=\ram(q_i), A_{ij}=\Irr(\Hom(\<q_i\>,\<q_j\>))$,

$\bullet$ If $i=j$ the number of arrows  from 
$\<q_i\>$ to $\<q_i\>$  (oriented loops) is 
\beq\label{eq: Bii} B_{ii}=A_{ii}-\be_i^2+1.\eeq

\begin{defn} 
The diagram $\wh \Ga$ of $(\Th,\cC)$ 
is obtained 
by gluing the end node of the leg $L_i$ to the node $\<q_i\>$ 
of the core diagram $\Ga$, for $i=1,\ldots,m$.
\end{defn}

As usual (when drawing diagrams) a pair of oppositely oriented arrows 
is identified with a single unoriented edge.
In particular a pair of oriented loops is the same thing as a single (unoriented) edge loop.
One can show (e.g. using the symplectic results of \cite{twcv})
that the integer $A_{ii}-\be_i^2+1$ is always even, and so 
all the diagrams here only involve unoriented edges (it is clear that $B_{ij}=B_{ji}$).
It should be noted that we call this a diagram, and not a graph, since some of the edges may have {\em negative} multiplicity.
(The negative edges will be indicated by dashed lines in figures.)
The meaning of a negative edge is that one has more relations than linear maps---the diagram arises since there are 
explicit matrix presentations of the wild character varieties, many of which date back to Birkhoff (cf. the history discussed in \cite{tops}).
It is something of a surprise that there are (lots of)
perfectly good moduli spaces whose diagrams have negative edges.

The untwisted case considered in \cite{rsode} is the case where each 
$\be_i=1$, so that $\<q_i\>\to \partial$ is a trivial (degree one) cover.
In this case each $q_i$ can be identified with a polynomial
in $x$ with zero constant term (where $x$ is a coordinate on 
$\IA^1=\IP^1\setminus \infty$).
In this case $\End(\<q_i\>)=\<0\>$ so that  $A_{ii}=0$ and so there are no edge loops. 
Further $\Hom(\<q_i\>,\<q_j\>) = \<q_j-q_i\>$
is again unramified, and its irregularity is just the degree of the 
polynomial $q_j-q_i$, so that 
\beq
A_{ij}-\be_i\be_j = \deg(q_i-q_j)-1.
\eeq
It follows that the core diagram coincides with the graph $\Ga$
defined 
in \cite{rsode} Apx C. The irregular type 
$\sum A_idz/z^{k-i}$ considered there 
is $dQ$ where
$Q$ is the diagonal matrix with entries given by the $q_i$, written in the coordinate $z=1/x$.
The simple expression $\deg(q_i-q_j)-1$ for the edge multiplicities of this graph %
appears in \cite{hi-ya-nslcase} \S3.3.

\subsection{Adding some tame singularities}\label{ssn: tame poles}

Here is how to extend this construction to the case where a finite number of tame singularities on $\IA^1$ are included as well. 
(This is similar to the procedure for adding tame singularities at finite distance in \cite{rsode, slims, cmqv}).
As mentioned in the introduction, this case 
(and any other case on $\IP^1$) can be reduced to 
the case already considered above via Fourier--Laplace.

Let $n=\rk(E)=\rk(V^0)=\sum_1^m n_i\be_i$ 
be the rank of the irregular class $\Th$ considered above. 
Choose points $a_1,\ldots,a_k\in \IA^1(\IC)$, and fix a tame 
formal class at each point. 
This is the same as fixing 
conjugacy classes 
$\wh \cC_1,\ldots, \wh \cC_k\subset \GL_n(\IC)$,
i.e. the local monodromy conjugacy classes.
(This is the same as fixing the isomorphism class of a graded local system at each point, but graded entirely by the corresponding 
tame circle $\<0\>$ with multiplicity $n$).

Let $\wh L_i$ be the leg determined by a 
minimal marking of 
$\wh \cC_i$ (in the sense of \cite{cmqv} Defn 9.2).
Assuming each conjugacy class is non-central,
each $\wh L_i$ has at least two nodes.

Let $\be= \sum \be_i=\ram(\Th)$ so that $\be \le n$.
Now  splay the end node of $\wh L_i$ into $\be$ nodes,  
thus replacing the end node by 
$\be$ nodes
(cf. \cite{rsode} Figure 6 \S A.5).
Glue the first $\be_1$ such nodes to the core node $\<q_1\>$.
Then glue the next $\be_2$ such nodes to the core node $\<q_2\>$, etc, thus gluing each of the $\be$ nodes to one of the core nodes.
In effect the second node of $\wh L_i$ is now linked to $\<q_j\>$ by 
$\be_j$ (unoriented) edges for $j=1,\ldots,m$.
Repeat this process for each $i=1,\ldots, k$.

This defines directly 
the diagram $\wh \Ga$ of any meromorphic connection $(E,\nabla)$ 
on $\IP^1$ that is tame at all but one point (i.e. associated to the choice of the formal class at $\infty$ plus tame 
classes at each $a_i$).

\section{Cartan matrix and dimensions}

Given a diagram $\wh \Ga$ with nodes $N$ and ``adjacency matrix'' $B$ 
(so $B_{ij}\in \IZ$ is the possibly negative number of arrows 
from node $i$ to $j$),
define the Cartan matrix of $\wh \Ga$ to be $C=2.\Id-B$.
Let $(\cdot,\cdot):\IZ^N\times \IZ^N\to \IZ$ be the resulting 
symmetric bilinear form defined by $(\eps_i,\eps_j)=C_{ij}$,
where the $\eps_i$ are the basis vectors of $\IZ^N$.
A {\em dimension vector}  for $\wh \Ga$ is a vector $\bd \in \IN^N$, i.e. a map $N\to \IN$  assigning a nonnegative integer to each node.

Now suppose $\wh \Ga$ is the diagram  
determined by a wild Riemann surface 
$\bSi=(\IP^1,\infty,\Th)$ 
with marked formal monodromy classes
$\cC = (\cC_1,\ldots,\cC_m)$ as above.
Then $\wh \Ga$ comes equipped with a dimension vector $\bd$
(the dimension of a core node $\<q_i\>$ is just 
$n_i=\Th(q_i)$ and the procedure described in
\cite{cmqv} Defn 9.2 gives the dimensions down the legs).
On the other  hand one can consider the  
(symplectic)  
wild character variety 
$$\MB(\bSi,\cC) \subset  \Hom_\IS(\Pi,G)/H$$
determined by this data, as in \cite{twcv} 
(except here we only consider the stable points).
It is the symplectic leaf of the Poisson wild character variety 
$\Hom_\IS(\Pi,G)/H$ determined by the classes $\cC$
(see also \cite{tops} which focuses on the general linear case).
Here $\Pi=\pi_1(\wt \Si,b)$ is the wild surface group 
(the 
fundamental group of the auxiliary surface $\wt \Si$), 
 $G=\GL_n(\IC)$ and $H=\Prod_1^m\GL_{n_i}(\IC)^{\be_i}\subset G$.
The classes $\cC$ determine a twisted conjugacy class of $H$ (by saying the monodromy around the circle $\<q_i\>$ is in $\cC_i$), and thus a symplectic leaf of $\Hom_\IS(\Pi,G)/H$.
(Fixing a symplectic leaf is the same as fixing the isomorphism class of the corresponding formal/$\cI$-graded local system $V^0$.)
The explicit presentations of the wild character varieties 
then leads to the statement:
\begin{prop}\label{prop:}
If $\MB(\bSi,\cC)$ is nonempty then 
$\dim_\IC(\MB(\bSi,\cC)) = 2-(\bd,\bd)$.
\end{prop}

The idea is that $\MB(\bSi,\cC)$ is thus a type of
multiplicative quiver variety for the ``doubled quiver $\wh \Ga$''.
This proposition can be proved directly, as sketched in the appendix.

\section{Examples}

It is easy to compute many examples.
Here are some of the simplest.
Note that if the multiplicities $n_i=1$ (so the formal monodromies are scalars) and there are no tame singularities, 
then we just need compute the core diagram.
Note also that everything is invariant under 
admissible deformations so we won't worry about constant factors---for example the (modern) Airy equation has  
$\Th=\<(2/3)x^{3/2}\>$, and the version $y''=9xy$ used  by Stokes in 1857 has 
$\Th=\<2x^{3/2}\>$, both of which are admissible deformations of 
$\<x^{3/2}\>$.
In the Painlev\'e cases below, unless explicitly stated otherwise, we use the linear equation in the standard Lax pair, due to Garnier/Jimbo--Miwa \cite{garn1912, JM81}.  (Recall the Painlev\'e equations arise as the isomonodromy equations of linear connections---the diagram of the Painlev\'e equation is that of the corresponding linear connection.)

$\bullet$\ Airy: $\Th=\<x^{3/2}\>$. 
$\End(\Th)=\<2x^{3/2}\>+2\<0\>$ so 
$A_{11}=3$ and in turn $B_{11} = 3-4+1=0$ so the Cartan matrix 
$C=(2)$ is the same as that for $\sl_2$, i.e. $A_1$.
The diagram has one node with zero loops (the $A_1$ Dynkin diagram).
$\dim(\MB)=0$.

$\bullet$\ Painlev\'e one: $\Th=\<x^{5/2}\>$.
$\End(\Th)=\<2x^{5/2}\>+2\<0\>$ so 
$A_{11}=5$ and in turn $B_{11} = 5-4+1=2$ so the Cartan matrix 
$C=(0)$ is the zero $1\times 1$ matrix.
The diagram has one node with one loop (the affine $A_0$ Dynkin diagram $\wh A_0$). $\dim(\MB)=2$.
This diagram appears also on the additive side---the corresponding additive moduli space $\cM^*$ is isomorphic to the affine plane 
$\IA^2$ (the $\wh A_0$ ALE space, familiar from the ADHM construction).

$\bullet$\ Weber: $\Th=\<x^2\>+\<-x^2\>$. This is untwisted so fits in to the set-up of \cite{rsode}.   
$A=\bsmx 0 & 2 \\ 2 & 0 \esmx$ so that $B=\bsmx 0 & 1 \\ 1 & 0 \esmx$
and $C=\bsmx 2 & -1 \\ -1 & 2 \esmx$. 
This is the Cartan matrix of $\sl_3$, i.e. $A_2$.
The diagram has two  nodes connected by one edge (the $A_2$ Dynkin diagram).
$\dim(\MB)=0$.

$\bullet$\ Painlev\'e two: $\Th=\<x^3\>+\<-x^3\>$.
This is untwisted so fits in to the set-up of \cite{rsode}.   
$A=\bsmx 0 & 3 \\ 3 & 0 \esmx$ so that $B=\bsmx 0 & 2 \\ 2 & 0 \esmx$
and $C=\bsmx 2 & -2 \\ -2 & 2 \esmx$. 
The diagram has two  nodes connected by two edges (the $\wh A_1$ Dynkin diagram).
This diagram appears also on the additive side 
(\cite{quad} Ex. 3)---the corresponding additive moduli space $\cM^*$ is diffeomorphic to the Eguchi--Hanson space (the $\wh A_1$ ALE space). $\dim(\MB)=2$.

$\bullet$\ Painlev\'e two revisited: $\Th=\<x^{3/2}\>$ plus a 
tame pole at $x=0$  (this is the Flaschka--Newell Lax pair, 
from the modified KdV equation).
The procedure of \S\ref{ssn: tame poles} again gives the $\wh A_1$ diagram: as in the Airy equation we get one node with no loops at $\infty$ (but with ramification $\be=2$). At the simple pole we get a leg of length $2$. We splay its end node into two nodes, and glue both of them to the node from $\infty$ yielding $\wh A_1$.

$\bullet$\ Bessel--Clifford equation (\!\!$\ _0F_1$-equation/confluent hypergeometric limit equation/Kummer's second equation, $xy'' + ay' =y$): 
$\Th=\<x^{1/2}\>$ plus a tame pole at $x=0$.
The procedure of \S\ref{ssn: tame poles} gives a diagram with two nodes attached by two edges. One node has no loops and the other has a single {\em negative} loop.
$C=\bsmx 2 & -2 \\ -2 & 4 \esmx$. $\dim(\MB)=0$.

$\bullet$\ Painlev\'e three: $\Th=\<x^{1/2}\>$ plus two 
tame poles at $x=0,1$. (This is the Lax pair for $P_3$ known as  
``degenerate Painlev\'e five'', \cite{JKT07} (6.17)).
The procedure of \S\ref{ssn: tame poles} gives a diagram with three nodes: two nodes each attached with two edges to a central node, and the central node has a single {\em negative} loop:

\begin{figure}[ht]
    \centering
\includegraphics[width=0.4\textwidth]{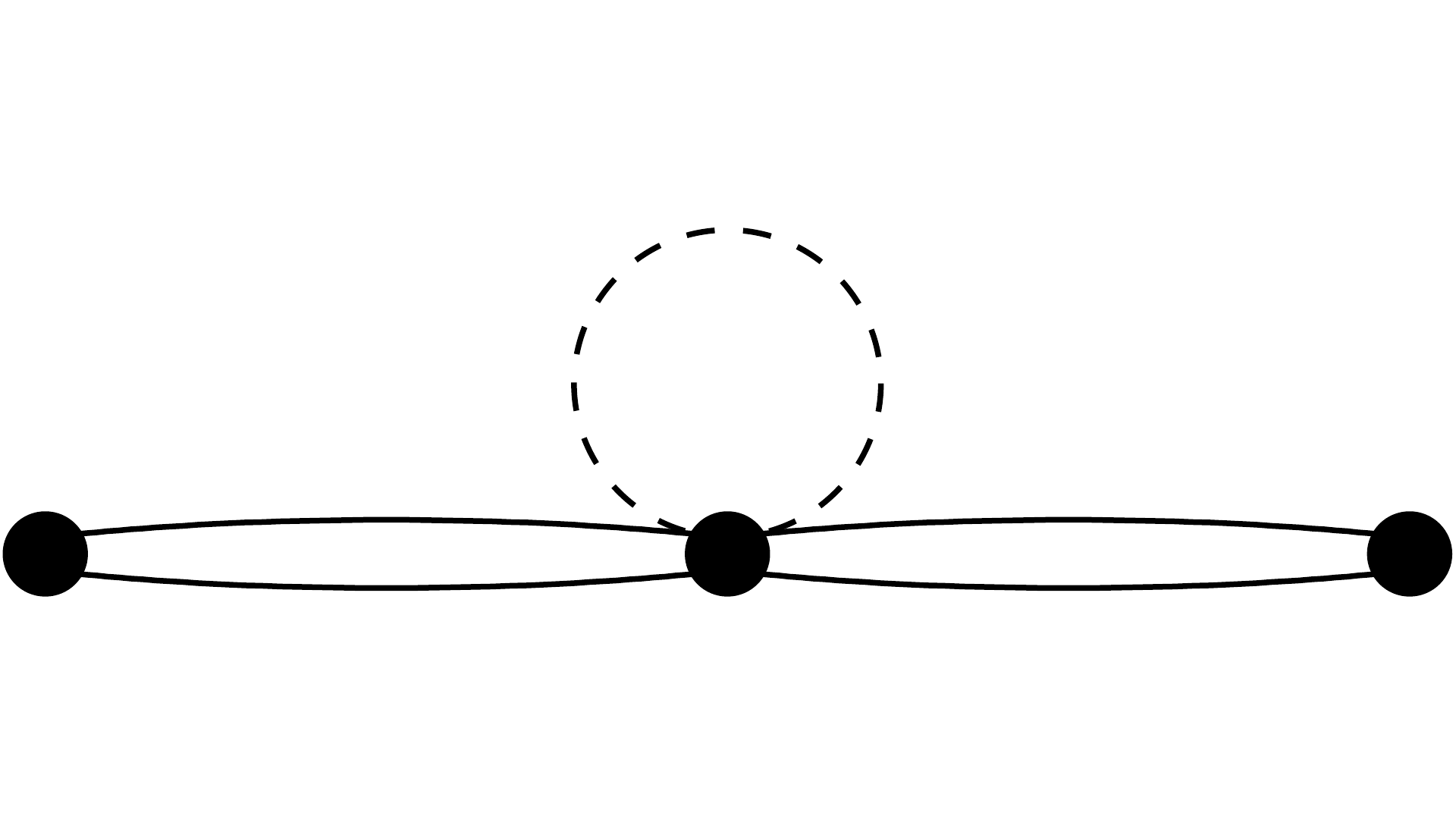}
\caption{The diagram of the Painlev\'e three equation.}
\end{figure}

The dashed line indicates that the loop has negative multiplicity.
The corresponding Cartan matrix is
\beq \label{eq: p3Cartan}
C=\bmx 2 & -2 & 0 \\ -2 & 4 & -2 \\ 0 & -2 & 2 \emx.
\eeq

The corresponding additive moduli space $\cM^*$ for Painlev\'e three 
(from the standard Lax pair) is known\footnote{This is stated in \cite{rsode}---it follows immediately since \cite{smid} implies $\cM^*\cong \IC^*\sqp T^*\SL_2(\IC)\spq \IC^*$, and this is 
the description of the 
$\wh D_2$ ALF space in \cite{dancer-dih} p.88.} 
to be the affine $D_2$ ALF space, and so it is natural to view this graph as the Dynkin diagram of type 
$\wh D_2$. 
As a further consistency check one can consider the 
intersection form from the corresponding De Rham moduli space 
(the Okamato space of initial conditions).  
It is known (\cite{Sakai-CMP01} p.182) that the intersection form 
is the negative of $\wh A_1\oplus A_1$, i.e. 
$-\bsmx  2 & -2 & 0 \\ -2 & 2 & 0 \\ 0 & 0 & 2 \\ \esmx$.
This fits since one can choose another $\IZ$-basis in which the intersection form is given by the negative of our  $\wh D_2$ Cartan matrix \eqref{eq: p3Cartan}:

\begin{lem}
The bilinear forms  with matrices
$\wh D_2$ and $\wh A_1\oplus A_1$ are equivalent 
over $\IZ$.
\end{lem}
\pf
$g^T\bsmx  2 & -2 & 0 \\ -2 & 2 & 0 \\ 0 & 0 & 2 \\ \esmx g = \bsmx 2 & -2 & 0 \\ -2 & 4 & -2 \\ 0 & -2 & 2 \esmx\qquad \text{where}\qquad
g=\bsmx  0 & 1 & 0 \\ 0 & 0 & 1 \\ -1 & 1 & 0 \\ \esmx.$
\epf

This diagram should also be compared/contrasted 
with the ``shape'' for Painlev\'e 3 suggested 
in \cite{hiroe-lde} Example 6.17 (and last shape in figure on p.928),
and with that in the approach of \cite{yamak-nonsym}.

The diagrams for Painlev\'e 4,5,6 have already been discussed in detail in \cite{quad, rsode, slims, cmqv}.
The diagrams for the six Painlev\'e equations are thus as follows:

\begin{figure}[ht]
    \centering
\includegraphics[width=0.9\textwidth]{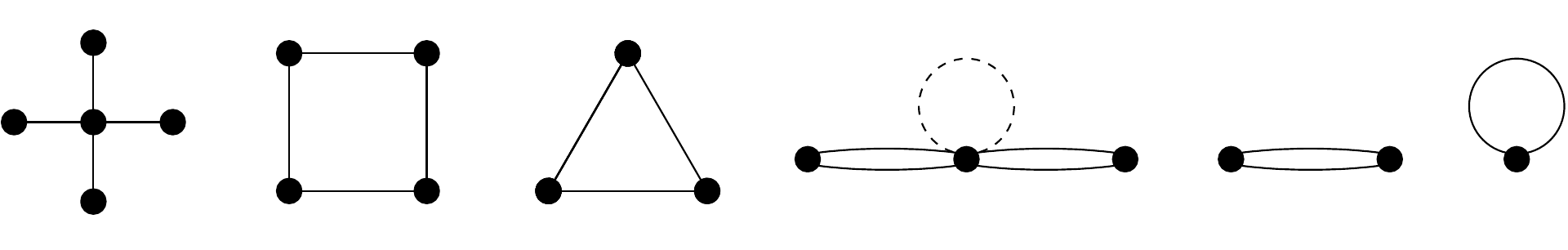}
\caption{The diagrams of the six Painlev\'e equations.}
\end{figure}

Note that the number of nodes minus one is always the number of parameters in the corresponding Painlev\'e equation.
Each node has dimension one, except for the central node in the 
Painlev\'e VI case, which has dimension $2$.

\begin{rmk}
Each of the Painlev\'e equations 2-6 has 
a special class of solutions 
coming from solutions of a linear differential 
equation (\cite{Okamoto-dynkin} p.310). 
We now understand the key statement:
{\em In each case the diagram of the linear equation is obtained by removing one node of the diagram for the corresponding Painlev\'e equation:} 

\noindent
  \begin{center}
  \begin{tabular}{| c || c | c | c | c | c | }
    \hline
Painlev\'e equation &  $P_6$ \phantom{$\Bigl\vert$}  &  $P_5$  &  $P_4$  &  $P_3$ & $P_2$\\ \hline 
Special solutions &  Gauss \!\!$\ _2F_1$\phantom{$\Bigl\vert$}    &  
Kummer \!\!$\ _1F_1$ &  Weber   &  
Bessel--Clifford \!\!$\ _0F_1$ &   Airy \\ \hline
  \end{tabular}
  \end{center}

Note that the discussion in \cite{rsode, slims} 
implies the Kummer equation has diagram $A_3$ and 
the Gauss equation has diagram $D_4$. 
Also \cite{Okamoto-dynkin} just writes ``Bessel'' for the special solutions of $P_3$, but the Bessel equation $x^2y'' + xy' =(\nu^2-x^2)y$ is really just a special one parameter
subfamily of the Kummer (\!\!$\ _1F_1$) equations 
(up to a twist)---it is the subfamily that are pullbacks of a \!\!$\ _0F_1$ equation (studied by Clifford):
If $f$ satisfies $xf'' + af' =f$ then 
$x^{a-1}\cdot f(-x^2/4)$ satisfies the Bessel 
equation with parameter $\nu=a-1$, for example:
$$
J_\nu(x) = \frac{(\frac{1}{2}x)^\nu}{\Ga(\nu+1)}.\ _0F_1(\nu+1; -x^2/4).$$
\end{rmk}

\appendix

\section{Sketch of proof of Prop. \ref{prop:}}

Recall from \cite{twcv} that the space
$\cB:=\Hom_\IS(\Pi,G)$ is isomorphic to 
$\cA\spq G$, where $\cA\cong G\times H(\partial)\times \ISto$ is the fission space, $H(\partial)$ is a twist of $H$ and $\ISto$ is the product of the Stokes groups (which has dimension $\Irr\End(\Th)$).
\cite{twcv} shows that $\cA$ is a twisted quasi-Hamiltonian $G\times H$ space, with a moment map $\mu=(\mu_G,\mu_H)$ 
to $G\times H(\partial)$. The $G$ action is free which implies 
$\dim \cB=\dim\cA-2\dim G$. 
In turn $\MB(\bSi,\cC)=\cB\spq_{\wt\cC} H=\mu_H^{-1}(\wt \cC)/H$ 
is the 
tq-Hamiltonian reduction (of the stable points)
by $H$ at the twisted conjugacy class 
$\wt \cC\subset H(\partial)$ determined by $\cC$.
This has dimension $\dim(\cB\times \wt \cC)-2\dim(H)+2$,
since $H/Z(G)$ acts effectively on stable points, and the result follows.

A better approach  is to frame the corresponding Stokes local systems slightly differently, as follows (this won't work for arbitrary reductive groups): Just choose one basepoint on each circle 
$\<q_i\>$, and frame there. The resulting space of framed Stokes local systems has the form $\cE:=\cB/H^\perp$ where 
$H^\perp\cong \Prod_1^m \GL_{n_i}(\IC)^{\be_i-1}\subset H$ (forgetting most of the old framings).
This has a residual action of 
$\check H:=\Prod_1^m \GL_{n_i}(\IC)$ (changing the remaining framings), and one can deduce from \cite{twcv} that $\cE$
is a quasi-Hamiltonian $\check H$-space, with moment map given by the formal monodromy, all the way around each circle $\<q_i\>$
(cf. \cite{tops} p.1---there is only one outer boundary circle).
Then $\MB(\bSi,\cC)=\cE\spq_{\cC} \check H$ is just the q-Hamiltonian reduction at the class $\cC\subset \check H$.
The space $\cE$ behaves as if it were the space of invertible representations of the ``doubled quiver'' given by the core $\Ga$---one can identify directly the positive 
terms in \eqref{eq: Bij},\eqref{eq: Bii}
with generators ({\em Stokes arrows} \cite{tops}, plus formal monodromies) and the negative terms with relations (from $\mu_G=1$). 
The dimension count is then standard, as in \cite{cmqv} \S9.1.

\noindent{\bf Acknowledgments.}
The second named author was supported by JSPS KAKENHI Grant Number 18K03256.

\renewcommand{\baselinestretch}{1}              %
\normalsize
\bibliographystyle{amsplain}    \label{biby}
\bibliography{../thesis/syr}

\vspace{0.5cm} 
\noindent
Laboratoire de Math\'ematiques d'Orsay, \\
Univ. Paris-Sud, CNRS, \\
Universit\'e Paris-Saclay, \\
91405 Orsay, France

\noindent
philip.boalch@math.u-psud.fr

\vspace{0.5cm} 
\noindent
Department of Mathematics,\\
Faculty of Science Division I,\\
Tokyo University of Science\\
1-3 Kagurazaka, Shinjuku-ku, \\
Tokyo 162-8601, Japan

\noindent
yamakawa@rs.tus.ac.jp

\end{document}